\documentclass{article}
\usepackage{graphicx}
\usepackage{amsfonts}
\usepackage{amsmath}
\usepackage{amssymb}
\usepackage{color}
\usepackage{algorithm}
\usepackage{algorithmic}
\usepackage{rotating}
\usepackage{textcomp}
\newtheorem{thm}{Theorem}[section]
\newtheorem{cor}{Corollary}[section]
\newtheorem{lem}[cor]{Lemma}

\newcommand{\R}{\mathbb R}

\newcommand{\N}{\mathbb N}

\newcommand{\A}{\mathcal{A}}

\newcommand{\OO}{\mathcal{O}}

\renewcommand{\L}{\mathcal{L}}

\newcommand{\la}{\left\langle}
\newcommand{\ra}{\right\rangle}
\newcommand{\lno}{\left\|}
\newcommand{\rno}{\right\|}
\newcommand{\sumweights}{\Gamma}
\newcommand{\de}{d_{\Delta_m}}

\newcommand{\dem}{d_{\Delta_n^M}}

\newcommand{\dgf}{distance-generating function}
\newcommand{\prox}{\textup{Prox}}

\newcommand{\qed}{${}$\hfill \rule{2mm}{2mm}}

\begin{document}
\title{An acceleration procedure for optimal first-order methods}
\author{Michel Baes\footnote{Corr. author. Institute for Operations Research, ETH Z\"urich, R\"amistrasse 101, 8092 Z\"urich, Switzerland,\newline michel.baes@ifor.math.ethz.ch.} , Michael B\"urgisser\footnote{Institute for Operations Research, ETH Z\"urich, R\"amistrasse 101, 8092 Z\"urich, Switzerland,\newline michael.buergisser@ifor.math.ethz.ch.}}
\maketitle

\begin{abstract} 
 We introduce in this paper an optimal first-order method that allows an easy and cheap evaluation of the local Lipschitz constant of the objective's gradient. This constant must ideally be chosen at every iteration as small as possible, while serving in an indispensable upper bound for the value of the objective function. In the previously existing variants of optimal first-order methods, 
 this upper bound inequality was constructed from points computed during the current iteration. It was thus not possible to select the optimal value for this Lipschitz constant at the beginning of the iteration. 
 
 In our variant, the upper bound inequality is constructed from points available before the current iteration, offering us the possibility to set the Lipschitz constant to its optimal value at once. This procedure, even if efficient in practice, presents a higher worse-case complexity than standard optimal first-order methods. We propose an alternative strategy that retains the practical efficiency of this procedure, while having an optimal worse-case complexity. We show how our generic scheme can be adapted for smoothing techniques, and perform numerical experiments on large scale eigenvalue minimization problems. As compared with standard optimal first-order methods, our schemes allows us to divide computation times by two to three orders of magnitude for the largest problems we considered.
 
 \end{abstract}

\section{Introduction}

With a few notable exceptions \cite{Gonzio_billions}, first-order methods constitute the main family of algorithms able to deal with very large-scale convex optimization problems \cite{Pena_Poker, nesterov:coreDP10/2010, Richtarik:block_coordinate, nesterov:coreDP16/2012}. Among them, optimal first-order methods play a distinguished role: they are practically as cheap as a first-order method can be, with a complexity per iteration growing as a moderate polynomial of the problem's size, while the worst-case number of iterations they require is provably optimal for \emph{smooth} instances \cite{Nemirovski_Yudin_book, nesterov_first_accel}. Their scope of applicability is restricted to optimization problems with a differentiable convex objective function $f$, whose gradient is globally Lipschitz continuous. Nevertheless, Nesterov introduced a systematic procedure, applicable to many nonsmooth convex functions, for building a smooth approximation to which one can apply an optimal first-order minimization 
algorithm \
cite{nesterov:coreDP12/2003}. His construction can easily be specified to realize the optimal compromise 
between the smoothness of the substitute objective function and how accurately it approximates the original objective. Smoothing techniques extended dramatically the scope of optimal first-order methods, and many variants of the original scheme developed in \cite{nesterov_first_accel} have been studied since then (see e.g. \cite{Tseng:proximal,Fista,d'Aspremont:Approx_gradients,Monteiro_restart}).

Critically, optimal first-order methods need an estimation of the corresponding Lipschitz constant with respect to an appropriate norm. Originally, this bound is used to build an approximation of the epigraph of the objective function. The larger the bound, the worse this approximation is, and the more steps the method is likely to take. Some strategies have been proposed to re-actualize at every step this bound \cite{nesterov_composite,Candes:first-order}. These strategies are based on the fact that the Lipschitz constant is used at a particular iteration to satisfy a single inequality rather than as a global property. If this inequality is verified, it suffices to reduce this constant, redo the iteration with the new value, and recheck the inequality, until it is no longer satisfied. If the inequality is not verified, we simply multiply the constant by an appropriate value and re-perform the iteration as long as the 
inequality does not hold. This 
strategy yielded a significant increase in practical efficiency. However, the cost of a single iteration has to be multiplied by a number ranging between two and possibly a few dozen. 

We show in this paper how a slight modification of these methods allows us to choose inexpensively the smallest possible approximation that guarantees the global convergence of the method. In particular, we avoid redoing several times the work needed for one iteration. The practical effect of such a procedure is appreciable, and is documented at the end of this paper. On the theoretical side, we show that our re-evaluation of the Lipschitz constant, if applied systematically, gives an algorithm which requires at worse $\OO((LD)/\epsilon)$ iterations, where $L$ is the global Lipschitz of the objective's gradient, $D$ measures the diameter of the feasible set, and $\epsilon>0$ is the desired absolute accuracy on the objective's value. In comparison with the vanilla optimal first-order method, which has a complexity of $\OO(\sqrt{(LD)/\epsilon})$ iterations, this algorithm is clearly worse. We propose a mixed strategy that presents simultaneously the practical efficiency of our systematic method for very large-
scale problems and, up to a constant that we can take as close to 1 as desired, the theoretical efficiency of optimal methods. 

When we apply to smoothing techniques, our mixed strategy suggests a different choice of the smoothness parameter than the standard one. This fact should not be too surprising: as our method is precisely designed to fit appropriate local estimates of the gradient's Lipschitz constant, it allows us to be slightly sloppier in our request for global smoothness. 

In order to validate our general scheme, we consider a well-known application of smoothing techniques to the problem of minimizing the largest eigenvalue of a convex combination of given symmetric matrices \cite{nesterov:coreDP73/2004}. This problem has many applications, and a large variety of methods have been devised to solve it \cite{Helmberg_bundle,Arora_Kale_07}, some of which are adaptations of optimal first-order methods \cite{Nem_Jud_Lan_Shap:stochastic,Buergisser_FastExp,d'Aspremont:stochastic_SDP}. To the best of our knowledge, these methods improve the complexity of each step of optimal first-order methods, but are not trying to decrease the number of these steps. With respect to original smoothing techniques, our method allows us to \emph{divide} by hundreds the practical number of iterations for large-scale instances, that is, when we have 100 matrices of size larger than 200. Our methods even allowed us to deal with a problem involving 10\% sparse matrices of dimension 12,800 within 9 hours, 
while standard optimal first-order methods would have taken more than one year if it were to perform all the iterations predicted by the worse-case analysis. It appears in practice that the standard optimal method needs about two third of these iterations: about eight months would be needed to solve that problem.

The paper is organized as follows. We outline our method in Section 2. First, we analyze its complexity for smooth convex problems and particularize our result to the two variants mentioned above. Then, we describe how the algorithm and its variants can be particularized to smoothed problems. In Section 3, we apply these methods to the eigenvalue minimization problem and present some numerical experiments. We have relegated the rather technical proof of the main theorem of the paper in the appendix.

\section{An accelerated optimal first-order method}
\label{sec:OFM}
In this section, we introduce an accelerated version of Nesterov's optimal first-order method that is presented in \cite{nesterov:coreDP12/2003} and discuss its application in smoothing techniques.

\subsection{General algorithm}

We start by considering the following optimization problem: 
\begin{eqnarray}\label{eq:opt_problem}f^*=\min_{x\in Q}f(x),\end{eqnarray}
where $Q$ is a closed and convex subset of $\R^n$ and $f:\R^n\rightarrow\R$ is a function, which is supposed to attain its minimum on the set $Q$. In addition, we assume that $f$ is convex and differentiable with a Lipschitz continuous gradient on $Q$. 

We consider $\R^n$ with the standard Euclidean scalar product, which is denoted by $\la\cdot,\cdot\ra$. The space $\R^n$ is equipped with a norm $\lno\cdot\rno_{\R^n}$, which may differ from the norm that is induced by the scalar product. We write $\lno\cdot\rno_{\R^n,*}$ for the dual norm to $\lno\cdot\rno_{\R^n}$:
\[\lno u\rno_{\R^n,*}:=\max_{x\in \R^n}\left\{ \la u, x\ra:\lno x\rno_{\R^n} =1\right\},\qquad u\in\R^n.\]

As $f$ has a Lipschitz continuous gradient on $Q$, there exists a constant $L=L(Q)>0$ which satisfies the inequality:
\begin{eqnarray}
 \label{eq:lipschitz_condition}
\lno \nabla f(x)-\nabla f(y)\rno_{\R^n,*}\leq L\lno x-y\rno_{\R^n}\qquad \forall\ x,y\in Q.
\end{eqnarray}
Nesterov developed a first-order method (see Equations (5.6) in \cite{nesterov:coreDP12/2003}) that allows us to compute approximate solutions to Problem (\ref{eq:opt_problem}). This optimal first-order method has a convergence rate of $\mathcal{O}\left(L/T^2\right)$, which outperforms the rate of convergence of common subgradient methods by two orders of magnitude. We quickly recall that common subgradient methods converge with the order $\OO(1/T^{0.5})$; see for instance \cite{Nemirovski_Yudin_book}. 

At every step of Nesterov's optimal first-order method, the Lipschitz constant $L$ is used to update the iterates; see \cite{nesterov:coreDP12/2003} for the details. However, the constant $L$ is a global parameter of the function $f$, as $L$ needs to satisfy Condition $(\ref{eq:lipschitz_condition})$ on the whole set $Q$. In this subsection, we introduce a refined version of Nesterov's optimal first-order method, where we replace the global parameter $L$ by local estimates. 

This algorithm requires the following basic notions. We say that $d_Q:Q\rightarrow \R_{\geq 0}$ is a \textit{\dgf\ for the set $Q$} if it complies with the following requirements:
\begin{enumerate}
 \item $d_Q$ is continuous on $Q$;
 \item $d_Q$ is strongly convex with modulus $1$ on $Q$:
 \[d_Q(\lambda x+[1-\lambda ]y)+\frac{\lambda [1-\lambda]}{2}\lno x- y\rno_{\R^n}^2\leq \lambda d_Q(x)+[1-\lambda]d_Q(y)\qquad\forall\ x,y\in Q;\]
\item given the set $Q^o(d_Q):=\left\{x \in Q:\partial d_Q(x)\neq \emptyset\right\}$, the subdifferential $\partial d_Q$ gives rise to a continuous selection $d'_Q$ on the set $Q^o$. If there is no possibility for confusion, we write $Q^o$ instead of $Q^o(d_Q)$.
\end{enumerate}

Let $d_Q$ be a \dgf\ for the set $Q$ and choose $z\in Q^o$. We write 
\[V_z^{d_Q}(x)=d_Q(x)-d_Q(z)-\left\langle d_Q'(z),x-z\right\rangle\in\R_{\geq 0}\] 
for the \textit{Bregman distance of $x\in Q$ with respect to $z\in Q^o$}. Nesterov's optimal first-order method and its accelerated version that we present in this paper utilize a \textit{prox-mapping}, that is, a mapping of the form:
\begin{eqnarray}
 \label{eq:prox_mapping}
\prox_{Q,z}^{d_Q}:\R^n\rightarrow Q^o:s\mapsto\arg\min_{x\in Q}\left\{\la s,x-z \ra+ V_z^{d_Q}(x)\right\},\qquad z\in Q^o.
\end{eqnarray}
If there is no possibility for confusion, we abbreviate $V_z^{d_Q}$ and $\prox_{Q,z}^{d_Q}$ into $V_z$ and  $\prox_{Q,z}$, respectively.
Given $s\in \R^n$ and $z\in Q^o$, the value $\prox_{Q,z}(s)$ can be rewritten as
\[\prox_{Q,z}(s)=\arg\min_{x\in Q}\left\{\la s-d_Q'(z),x \ra+ d_Q(x)\right\}.\]
It can be easily verified that this optimization problem has indeed a unique minimizer (Note that the objective function $x\mapsto\la s-d_Q'(z),x \ra+ d_Q(x)$ is continuous and strongly convex. It remains to apply Lemma $6$ from \cite{nesterov:PDS}.) and that this minimizer belongs to $Q^o$. For the reminder of this paper, we assume that this minimizer can be computed easily (Ideally, we can write it in a closed form.). The unique element 
\[c(d_Q):=\arg\min_{x\in Q}\left\{d_Q(x)\right\}\in Q^o\]
is called the $d_Q$-center (Note that $c(d_Q)=\prox_{Q,z}(d_Q'(z))$ for any $z\in Q^o$.). Without loss of generality, we may assume that $d_Q$ vanishes at the point $c(d_Q)$. Then,  Lemma 6 in \cite{nesterov:PDS} can be used to justify the following inequality:
\begin{eqnarray}\label{eq:lower_bound_dgf}
d_Q(x)\geq \frac{1}{2}\lno x-c(d_Q)\rno_{\R^n}^2\qquad\forall\ x\in Q. 
\end{eqnarray}

We discuss now the analytical complexity of the accelerated optimal first-order method displayed in Algorithm \ref{alg:opt_first_order}. We choose $T\in\N_0:=\N\cup\left\{ 0\right\}$  and assume that the sequences $\left(x_t\right)_{t=0}^{T+1}$, $\left(u_t\right)_{t=0}^{T+1}$, $\left(z_t\right)_{t=0}^T$, $\left(\hat{x}_t\right)_{t= 1}^{T+1}$, $\left(\gamma_t\right)_{t=0}^{T+1}$, $\left(\sumweights_t\right)_{t= 0}^{T+1}$, $\left(\tau_t\right)_{t=0}^{T}$, and $\left(L_t\right)_{t=0}^T$ are generated by Algorithm \ref{alg:opt_first_order}. Given $0\leq t\leq T$, we say that \textit{Inequality $(\mathcal{I}_t)$ holds} if 
\begin{equation}
\sumweights_tf(u_t)+\sum_{k=0}^{t-1}\left(L_{k+1}-L_k\right)\left(d_Q(z_{k+1})-\frac{1}{2}\left\|z_k-\hat{x}_{k+1}\right\|_{\R^n}^2 \right)\leq \psi_t, \tag{$\mathcal{I}_t$}
\end{equation}
where
\[\psi_t:=\min_{x\in Q}\left\lbrace
       \sum_{k=0}^{t}\gamma_k\left(f(x_k)+\langle \nabla f(x_k), x-x_k\rangle \right) + L_td_Q(x) \right\rbrace.\]

\begin{algorithm}[t]
  \caption{Accelerated optimal first-order method}
  \label{alg:opt_first_order}
  \begin{algorithmic}[1]
    \STATE Choose $T\in\N_0$.
    \STATE Choose $\left(\gamma_t\right)_{t=0}^{T+1}$ with $\gamma_0\in (0,1]$, $\gamma_t\geq 0$, and $\gamma_t^2\leq \sumweights_t:=\sum_{k=0}^t\gamma_k$ for any $0\leq t \leq T+1$.
    \STATE Set $L_0=L$ and $x_0=c(d_Q)$.
    \STATE Compute $u_0:=\arg\min_{x\in Q}\left\lbrace
   \gamma_0 \left(f(x_0)+\langle \nabla f(x_0), x-x_0\rangle \right)+  L_0d_Q(x)\right\rbrace$.
   \STATE Set $z_0=u_0$, $\tau_0=\gamma_1/\sumweights_1$, and $x_1=\tau_0z_0+(1-\tau_0)u_0=z_0$.
  \STATE Define $\hat{x}_1:=\prox_{Q,z}\left(\gamma_{1}\nabla f(x_1)/L_0\right)$.
\STATE Set $u_1=\tau_0\hat{x}_1+(1-\tau_0)u_0.$
    \FOR {$1\leq t \leq T$}
      \STATE Choose $0<L_t\leq L$ such that:
\begin{eqnarray}\label{eq:adapted_L_cond_opt}
      f(u_t)\leq f(x_t)+\left\langle \nabla f(x_t),u_t-x_t\right\rangle+\frac{L_t}{2}\left\| u_t-x_t\right\|_{\R^n}^2.\end{eqnarray}\vspace{-2mm}
\STATE Set $z_t=\arg\min_{x\in Q}\left\lbrace\sum_{k=0}^{t}\gamma_k\left(f(x_k)+\langle \nabla f(x_k), x-x_k\rangle \right)+L_td_Q(x) \right\rbrace$.\hspace{-1mm}
\STATE Set $\tau_t=\gamma_{t+1}/\sumweights_{t+1}$ and $x_{t+1}=\tau_t z_t+(1-\tau_t)u_t$.
\STATE Compute $\hat{x}_{t+1}:=\prox_{Q,z_t}\left(\gamma_{t+1}\nabla f(x_{t+1})/L_{t}\right)$.
        \STATE Set $u_{t+1}=\tau_t\hat{x}_{t+1}+(1-\tau_t)u_t.$
    \ENDFOR
    \end{algorithmic}
\end{algorithm}

As the proof of the following result is rather long and technical, we give it in the Appendix \ref{sec:proof_opt_first_order}.

\begin{thm}\label{thm:I_t}
 Inequality $(\mathcal{I}_t)$ holds for any $0\leq t\leq T$.
\end{thm}

For the reminder of this subsection, we refer to $x^*\in Q$ as an optimal solution to the optimization problem $f^*=\min_{x\in Q}f(x)$. 

\begin{thm}\label{thm:convergence_OFM}
For any $T\in\N_0$, we have:
\[f(u_T)-f^*\leq\frac{1}{\sumweights_T}\left[L_T d_Q(x^*)+\sum_{t=0}^{T-1}\left(L_t-L_{t+1}\right)\left(d_Q(z_{t+1})-\frac{1}{2}\left\|z_t-\hat{x}_{t+1}\right\|_{\R^n}^2  \right)\right].\]
\end{thm}

\noindent\textbf{Proof:} Let $0\leq t\leq T$. The convexity of the function $f$ and the definition of $\sumweights_t$ imply
\begin{eqnarray*}
 \psi_t&:=&\min_{x\in Q}\left\lbrace
        L_t d_Q(x)+\sum_{k=0}^{t}\gamma_k\left(f(x_k)+\langle \nabla f(x_k), x-x_k\rangle \right) \right\rbrace\cr
&\leq& L_td_Q(x^*)+\sum_{k=0}^{t}\gamma_k\left(f(x_k)+\langle \nabla f(x_k), x^*-x_k\rangle \right) \cr
&\leq& L_td_Q(x^*)+\sum_{k=0}^{t}\gamma_kf(x^*)\cr
&=&L_t d_Q(x^*)+\sumweights_t f(x^*).
\end{eqnarray*}
It remains to combine this inequality with Theorem \ref{thm:I_t}.
\qed

Nesterov \cite{nesterov:coreDP12/2003} suggests to choose the sequence $\left(\gamma_t\right)_{t=0}^{T+1}$ as
\begin{eqnarray}
 \label{eq:choice_gamma_t_s}
\gamma_t:=\frac{t+1}{2}\qquad \forall\ 0\leq t\leq T+1.
\end{eqnarray}
Lemma 2 of \cite{nesterov:coreDP12/2003} shows that we have the following equations for this choice of the sequence $\left(\gamma_t\right)_{t=0}^{T+1}$:
\[\tau_t=\frac{2}{t+3}\qquad\forall\ 0\leq t\leq T\]
and 
\[\sumweights_t=\frac{(t+1)(t+2)}{4},\qquad \gamma_t^2\leq \sumweights_t\qquad \forall\ 0\leq t\leq T+1.\]
As an immediate consequence of Theorem \ref{thm:convergence_OFM}, we obtain the following result for our accelerated optimal first-order method.
\begin{cor}
 Let us choose the sequence $(\gamma_t)_{t=0}^{T+1}$ in Algorithm \ref{alg:opt_first_order} as described in (\ref{eq:choice_gamma_t_s}). Then, we have for any $T\in\N_0$:
\begin{eqnarray}\label{eq:conv_ofom}
f(u_T)-f^*\leq \frac{4L_Td_Q(x^*)}{(T+1)(T+2)}+\sum_{t=0}^{T-1}\frac{4\left(L_t-L_{t+1}\right)}{(T+1)(T+2)}\left(d_Q(z_{t+1})-\frac{1}{2}\left\|z_t-\hat{x}_{t+1}\right\|_{\R^n}^2  \right).
\end{eqnarray}
\end{cor}
\qed

There exist different strategies for updating the sequence $\left(L_t\right)_{t=0}^T$ in Algorithm \ref{alg:opt_first_order}. When $\alpha:=0$, we recover the complexity results of Nesterov's optimal first-order method (see e.g. Subsection 5.3 of \cite{nesterov:coreDP12/2003}), for which Inequality (\ref{eq:conv_ofom}) can be rewritten as
$f(u_T)-f^*\leq\frac{4Ld_Q(x^*)}{(T+1)(T+2)}$.

\textbf{Alternative 1: (most aggressive adaptive setting)} Fix $0<\kappa\ll 1$ and let $1\leq t\leq T$. The most aggressive choice for the constant $L_t$ corresponds to 
\begin{eqnarray}
 \label{eq:L_most_aggressive}
    L_t:=\max\left\{\bar L_t, \kappa L\right\}\in [\kappa L,L],\qquad \bar L_t:= \frac{2\left[f(u_t)-f(x_t)-\left\langle \nabla f(x_t),u_t-x_t\right\rangle\right]}{\left\| u_t-x_t\right\|_{\R^n}^2}\leq L.
\end{eqnarray}
The computation of the constant $L_t$ requires the entities $u_t$, $x_t$, and $\nabla f(x_t)$. In sharp contrast with the methods proposed so far \cite{nesterov_composite,Candes:first-order}, all these entities are known from the previous step $t-1$, implying that the constant $L_t$ can be determined immediately.

Independent of the choice of the $L_t$'s, we can always derive the following trivial convergence result for Algorithm \ref{alg:opt_first_order} from Inequality (\ref{eq:conv_ofom}): 
\begin{eqnarray*}
f(u_T)-f^*\leq \frac{4L\sup_{x\in Q}d_Q(x)}{(T+1)(T+2)}+\frac{20 LT \sup_{x\in Q}d_Q(x)}{(T+1)(T+2)}\leq \frac{20 L \sup_{x\in Q}d_Q(x)}{T+2},
\end{eqnarray*}
as $L_Td(x^*)\leq L \sup_{x\in Q}d_Q(x)$ and
\begin{eqnarray*}
\sum_{t=0}^{T-1}\left(L_t-L_{t+1}\right)\left(d_Q(z_{t+1})-\frac{1}{2}\left\|z_t-\hat{x}_{t+1}\right\|_{\R^n}^2\right)&\leq& \sum_{t=0}^{T-1}\left|L_t-L_{t+1}\right|\left(d_Q(z_{t+1})+\frac{1}{2}\left\|z_t-\hat{x}_{t+1}\right\|_{\R^n}^2\right)\cr
&\leq& \sum_{t=0}^{T-1}L \left(\sup_{x\in Q}d_Q(x)+4\sup_{x\in Q}d_Q(x)\right)\cr
&=&5 L T \sup_{x\in Q}d_Q(x).
\end{eqnarray*}
Note that the last inequality holds due to (\ref{eq:lower_bound_dgf}).

Thus, Algorithm \ref{alg:opt_first_order} equipped with the most aggressive update strategy, which is described in (\ref{eq:L_most_aggressive}), needs at most \[T=\left\lceil 20 L\sup_{x\in Q}d(x)/\epsilon -2\right\rceil\] iterations to find a feasible $\epsilon$-solution, provided that $\sup_{x\in Q}d(x)$ is finite. 

\textbf{Alternative 2: (hybrid setting)} Finally, we can combine the two settings that are presented above. We choose a number $\alpha\geq 0$ and denote by $1\leq t\leq T$ the current iteration. As long as  
\begin{eqnarray}\label{eq:switch_back_rule} 
\sum_{k=0}^{\bar t-1}\left(L_k-L_{k+1}\right)\left(d_Q(z_{k+1})-\frac{1}{2}\left\|z_k-\hat{x}_{k+1}\right\|_{\R^n}^2  \right)\leq \alpha L d_Q(x^*)\qquad \forall\ 1\leq \bar t\leq t,
\end{eqnarray}
we use the update strategy that is described in (\ref{eq:L_most_aggressive}). When Condition (\ref{eq:switch_back_rule}) is not satisfied for the first time, we set $L_{\bar t}:=L$ for any $\bar t\geq t$ and recompute the point $z_t$. 

With the just specified setting, Inequality (\ref{eq:conv_ofom}) results in the bound
\[f(u_T)-f^*\leq \frac{4(1+\alpha) L d_Q(x^*)}{(T+1)(T+2)}.\]
That is, we need to perform at most 
\[T=\left\lceil 2\sqrt{(1+\alpha) Ld_Q(x^*)/\epsilon}-1\right\rceil\]
iterations of Algorithm \ref{alg:opt_first_order} to find a point $x\in Q$ with $f(x)-f^*\leq \epsilon$, where $\epsilon >0$. This complexity result deviates by a factor of $(1+\alpha)^{0.5}$ from the efficiency estimate of the non-adaptive method. With $\alpha=5(T+1)\sup_{x\in Q}d(x)-1$, the setting coincides with Alternative 1.

\subsection{The accelerated optimal first-order method in smoothing techniques}
\label{sec:ST}

Smoothing techniques \cite{nesterov:coreDP12/2003} constitute a two-stage procedure that can be applied to non-smooth optimization problems with a very particular structure. In a first step, a smooth approximation of the non-smooth objective function is formed, so that Nesterov's optimal first-order method can be applied afterwards. In this section, we study the effects of replacing Nesterov's original optimal first-order method by its accelerated version in smoothing techniques.

We assume that the sets $Q_1\subset\R^n$ and $Q_2\subset\R^m$ are both compact and convex. In addition, we endow the spaces $\R^n$ and $\R^m$ with two (maybe different) norms. We denote by $\lno\cdot\rno_{\R^n}$ and $\lno\cdot\rno_{\R^m}$ the norm of the spaces $\R^n$ and $\R^m$, respectively. Nesterov considers convex optimization problems of the form:
\begin{eqnarray}
 \label{eq:problem_min_max}
\min_{x\in Q_1}\max_{y\in Q_2} \phi(x,y),\qquad \phi(x,y):= f_1(x) + \la \mathcal{A}(x),y\ra - f_2(y),
\end{eqnarray}
where $f_1:\R^n\rightarrow \R$, $f_2:\R^m\rightarrow\R$ are smooth and convex, and $\mathcal{A}:\R^n\rightarrow\R^m$ is a linear operator. With a slight abuse of notation, we write $\la\cdot,\cdot\ra$ for the Euclidean scalar product in both spaces $\R^n$ and $\R^m$. 

According to the standard MiniMax Theorem in Convex Analysis (see Corollary 37.3.2 in \cite{Rockafellar:70}), we have, due to the compactness and convexity of the sets $Q_1$ and $Q_2$, the following pair of primal-dual convex optimization problems:
\[\min_{x\in Q_1}\left\{\overline{\phi}(x):=\max_{y\in Q_2}\phi(x,y)\right\}=\max_{y\in Q_2}\left\{\underline{\phi}(y):=\min_{x\in Q_1}\phi(x,y)\right\}.\]
The operator $\mathcal{A}$ comes with an adjoint operator $\mathcal{A}^*:\R^m\rightarrow\R^n$, which is defined by the relation:
\[\la \mathcal{A}(x),y\ra = \la x,\mathcal{A}^*(y)\ra\qquad \forall\ (x,y)\in \R^n\times\R^m.\]
The analysis of Nesterov's smoothing techniques requires a norm of the operator $\mathcal{A}$. This norm is constructed as follows: 
\[\lno \mathcal{A}\rno_{\R^n,\R^m}:=\max_{x\in\R^n,y\in\R^m}\left\{\la\A(x),y\ra:\lno x\rno_{\R^n}=1,\ \lno y\rno_{\R^m}=1\right\}.\] 

We are ready to form a smooth approximation of $\overline{\phi}$ to which we can apply Algorithm \ref{alg:opt_first_order}. We choose a \dgf\ $d_{Q_2}:Q_2\rightarrow\R_{\geq 0}$ for the set $Q_2$ and consider the auxiliary function
\[\overline{\phi}_\mu:\R^n\rightarrow\R:x\mapsto \max_{y\in Q_2}\left\{f_1(x)+\la \mathcal{A}(x),y\ra-f_2(y)-\mu d_{Q_2}(y)\right\},\]
where $\mu>0$ is a positive smoothness parameter. This function defines a uniform approximation of $\overline{\phi}$, as 
\begin{eqnarray}\label{eq:smooth_approximation_bounds}
\overline{\phi}_\mu(x)\leq \overline{\phi}(x)\leq \overline{\phi}_\mu(x)+\mu \max_{z\in Q_2}d_{Q_2}(z)\qquad \forall\ x\in Q_1; 
\end{eqnarray}
see Inequality (2.7) in \cite{nesterov:coreDP12/2003}. The function $y\mapsto \la \mathcal{A}(x),y\ra-f_2(y)-\mu d_{Q_2}(y)$ is strongly concave for any $x\in Q_1$, as the \dgf\ $d_{Q_2}$ is strongly convex by its definition. Hence, the function $y\mapsto \la \mathcal{A}(x),y\ra-f_2(y)-\mu d_{Q_2}(y)$ has a unique maximizer on $Q_2$. We denote this maximizer by $y_*(x)$. 

Nesterov showed that $\overline{\phi}_\mu$ is differentiable with a Lipschitz continuous gradient. We write $M>0$ for the Lipschitz constant of the gradient of $f_1$.

\begin{thm}[Theorem 1 in \cite{nesterov:coreDP12/2003}] \label{thm:gradient} The function $\overline{\phi}_\mu$ is well-defined, continuously differentiable, and convex on $\R^n$. The gradient of $\overline{\phi}_\mu$ takes the form 
\[\nabla \overline{\phi}_\mu(x)=\nabla f_1(x)+\A^*(y_*(x)),\]
and is Lipschitz continuous with the constant $L_\mu:=M+\lno \A\rno_{\R^n,\R^m}^2/\mu$. 
\end{thm}

As an immediate consequence, we can apply Algorithm \ref{alg:opt_first_order} to the problem:
\begin{eqnarray}
\label{eq:smooth_problem}
\min_{x\in Q_1}\overline{\phi}_\mu(x). 
\end{eqnarray}

\begin{algorithm}[t]
  \caption{Algorithm \ref{alg:opt_first_order} (with $\gamma_t=(t+1)/2$) applied to Problem (\ref{eq:smooth_problem})}
  \label{alg:smoothing_techniques}
  \begin{algorithmic}[1]
    \STATE Choose $T\in\N_0$.
    \STATE Choose a smoothness parameter $\mu>0$ and a \dgf\ $d_{Q_1}:Q_1\rightarrow\R$ for the set $Q_1$.
    \STATE Set $L_0=L_\mu=M+\lno \A\rno_{\R^n,\R^m}^2/\mu$ and $x_0=c(d_{Q_1})$.
    \STATE Set $u_0=\arg\min_{x\in Q_1}\left\lbrace
   \frac{1}{2} \left(\overline{\phi}_\mu(x_0)+\langle \nabla \overline{\phi}_\mu(x_0), x-x_0\rangle \right)+  L_0d_{Q_1}(x)\right\rbrace$.
   \STATE Set $z_0=u_0$, $\tau_0=\frac{2}{3}$, and $x_1=\tau_0z_0+(1-\tau_0)u_0=z_0$.
  \STATE Define $\hat{x}_1:=\prox_{Q_1,z}\left(\nabla \overline{\phi}_\mu(x_1)/L_0\right)$.
\STATE Set $u_1=\tau_0\hat{x}_1+(1-\tau_0)u_0.$
    \FOR {$1\leq t\leq T$}
      \STATE Choose $0<L_t\leq L_\mu$ such that:
\begin{eqnarray*}
	  \overline{\phi}_\mu(u_t)\leq \overline{\phi}_\mu(x_t)+\left\langle \nabla \overline{\phi}_\mu(x_t),u_t-x_t\right\rangle+\frac{L_t}{2}\left\| u_t-x_t\right\|_{\R^n}^2.
\end{eqnarray*}\vspace{-2mm}
        \STATE Set \[z_t=\arg\min_{x\in Q_1}\left\lbrace
        \sum_{k=0}^{t}\frac{k+1}{2}\left(\overline{\phi}_\mu(x_k)+\langle \nabla \overline{\phi}_\mu(x_k), x-x_k\rangle \right)+L_td_{Q_1}(x) \right\rbrace.\]\vspace{-2mm}
        \STATE Set $\tau_t=\frac{2}{t+3}$ and $x_{t+1}=\tau_t z_t+(1-\tau_t)u_t$.
        \STATE Compute $\hat{x}_{t+1}:=\prox_{Q_1,z_t}\left(\frac{t+2}{2}\nabla \overline{\phi}_\mu(x_{t+1})/L_{t}\right)$.
        \STATE Set $u_{t+1}=\tau_t\hat{x}_{t+1}+(1-\tau_t)u_t.$
    \ENDFOR
    \end{algorithmic}
\end{algorithm}

Algorithm \ref{alg:smoothing_techniques} corresponds to Algorithm \ref{alg:opt_first_order} when we apply this method with step-sizes as described in (\ref{eq:choice_gamma_t_s}) to Problem (\ref{eq:smooth_problem}). A slight adaptation of the proof of Theorem 3 in \cite{nesterov:coreDP12/2003} yields to the following result, for which we need the definitions:
\[D_1:=\max_{x\in Q_1}d_{Q_1}(x)\qquad \text{and}\qquad D_2:=\max_{y\in Q_2}d_{Q_2}(y).\]

\begin{thm}\label{thm:convergence_smoothing}
 Fix $T\in\N_0$ and assume that the sequences $\left(x_t\right)_{t= 0}^{T+1}$, $\left(u_t\right)_{t= 0}^{T+1}$, $\left(z_t\right)_{t= 0}^T$, $\left(\hat{x}_t\right)_{t= 1}^{T+1}$, and $\left(L_t\right)_{t= 0}^T$ are generated by Algorithm \ref{alg:smoothing_techniques} with the smoothness parameter $\mu>0$. For
\[\bar{x}:=u_T\in Q_1\qquad\text{and}\qquad \bar{y}:=\sum_{t=0}^T\frac{2(t+1)}{(T+1)(T+2)}y_*(x_t)\in Q_2,\]
we have:
\begin{eqnarray}\label{eq:bound_st} \overline{\phi}(\bar{x})-\underline{\phi}(\bar{y})\leq \frac{4\left(D_1\lno \A\rno_{\R^n,\R^m}^2/\mu+D_1M-\chi_T\right) }{(T+1)^2}+\mu D_2,\end{eqnarray}
where 
\[\chi_T:=\sum_{t=0}^{T-1}\left(L_{t+1}-L_t\right)\left(d_{Q_1}(z_{t+1})-\frac{1}{2}\left\|z_t-\hat{x}_{t+1}\right\|_{\R^n}^2 \right).\]
\end{thm}

For remainder of this section, we use the notations of Algorithm \ref{alg:opt_first_order} and Theorem \ref{thm:convergence_smoothing}.

\textbf{Proof:} In accordance to Theorem \ref{thm:I_t} and to the step-size choice (\ref{eq:choice_gamma_t_s}), we have the inequality:
\begin{eqnarray}
\label{eq:intermediate_smoothing_techniques_1}
 \overline{\phi}_\mu(\bar{x})= \overline{\phi}_\mu(u_T)
\leq \frac{4(L_TD_1-\chi_T)}{(T+1)(T+2)}+\min_{x\in Q_1}\frac{2\beta_T(x)}{(T+1)(T+2)},
\end{eqnarray}
where 
\[\beta_T(x):= \sum_{t=0}^T(t+1)\left(\overline{\phi}_\mu(x_t)+\la \nabla \overline{\phi}_\mu(x_t),x-x_t\ra\right)\qquad\forall\ x\in Q_1.\]
Let $x\in Q_1$. Using Theorem \ref{thm:gradient} and the convexity of $f_1$ and $f_2$, we can write:
\begin{eqnarray*}
\beta_T(x)&= &\sum_{t=0}^T(t+1)\left(f_1(x)+\la \A(x),y_*(x_t)\ra-f_2(y_*(x_t))-\mu d_{Q_2}(y_*(x_t)) \right)\cr
&\leq& \sum_{t=0}^T(t+1)\left(f_1(x)-\la \A(x),y_*(x_t)\ra -f_2(y_*(x_t))\right)\cr
&\leq&\frac{(T+1)(T+2)}{2}\left( f_1(x) +\la \A(x),\bar{y}\ra -f(\bar y)\right).
\end{eqnarray*}
The above inequality implies:
\begin{eqnarray}
\label{eq:intermediate_smoothing_techniques_2}
\min_{x\in Q_1}\beta_T(x)\leq \frac{(T+1)(T+2)}{2} \underline{\phi}(\bar{y}). 
\end{eqnarray}
Recall that we have $L_T\leq L_\mu=\lno\mathcal{A}\rno_{\R^n,\R^m}^2/\mu+M$ by construction. We use Inequalities (\ref{eq:intermediate_smoothing_techniques_1}), (\ref{eq:intermediate_smoothing_techniques_2}), and (\ref{eq:smooth_approximation_bounds}) to justify the following inequalities:
\begin{eqnarray*}\label{eq:aux}
 \frac{4(D_1\lno\A\rno_{\R^n,\R^m}^2/\mu-\chi_T)}{(T+1)^2}\geq \frac{4(L_TD_1-\chi_T)}{(T+1)(T+2)}
\geq\overline{\phi}_\mu(\bar{x})-\underline{\phi}(\bar{y})\geq\overline{\phi}(\bar{x})-\underline{\phi}(\bar{y})-\mu D_2.
\end{eqnarray*}
\qed

We conclude this section by discussing different strategies for choosing the sequence $(L_t)_{t=0}^T$ and the smoothness parameter $\mu$.

\textbf{Alternative 1: (most aggressive adaptive setting)} We can always give the following upper bound for the quantity $(-\chi_T)$ in Theorem \ref{thm:convergence_smoothing}:
\[-\chi_T\leq  5 L_\mu D_1 T =5 D_1T\left(\lno\A\rno_{\R^n,\R^m}^2/\mu+M\right),\]
which allows us to reformulate (\ref{eq:bound_st}) as
\[\overline{\phi}(\bar{x})-\underline{\phi}(\bar{y})\leq \frac{20D_1\left(\lno\A\rno_{\R^n,\R^m}^2/\mu +M\right)}{T+1}+\mu D_2.\]
Minimizing the right-hand side of the above inequality with respect to $\mu$, that is, setting $\mu$ to 
\[\mu_2^*:=2\lno \A\rno_{\R^n,\R^m}\sqrt{\frac{5 D_1}{(T+1)D_2}},\]
we obtain:
\[\overline{\phi}(\bar{x})-\underline{\phi}(\bar{y})\leq 4 \lno \A\rno_{\R^n,\R^m}\sqrt{\frac{5 D_1D_2}{T+1}}+\frac{20D_1M}{T+1}.\]
As this bound is independent of the choice the $L_t$'s, it is valid also for the most aggressive setting, that is, for 
\begin{eqnarray}
 \label{eq:L_most_aggressive_st}
    L_t:=\max\left\{\bar L_t, \kappa L_\mu\right\}\in [\kappa L_\mu,L_\mu],\qquad \bar L_t:= \frac{2\left[\overline{\phi}(u_t)-\overline{\phi}(x_t)-\left\langle \nabla \overline{\phi}(x_t),u_t-x_t\right\rangle\right]}{\left\| u_t-x_t\right\|_{\R^n}^2}\leq L_\mu,
\end{eqnarray}
where $1\leq t\leq T$ and $0<\kappa\ll 1$ is fixed.

\textbf{Alternative 2: (hybrid setting)} Let $\alpha\geq 0$. We follow the setting described in (\ref{eq:L_most_aggressive_st}) for all $1\leq t\leq T$ as long as
\[-\chi_{\bar t}\leq \alpha D_1\left(\lno \A\rno_{\R^n,\R^m}^2/\mu+M\right)\]
is satisfied for any $1\leq\bar t\leq t$. When this condition fails for the first time, say for $t=t'$, we set $L_t$ to $L_\mu$ for any $t\geq t'$ and recompute the point $z_{t'}$. In this hybrid setting, Inequality (\ref{eq:bound_st}) yields to 
\[\overline{\phi}(\bar{x})-\underline{\phi}(\bar{y})\leq \frac{4(1+\alpha)D_1 \left(\lno\A\rno_{\R^n,\R^m}^2/\mu+M\right)}{(T+1)^2}+\mu D_2.\]
We choose $\mu$ such that the right-hand side of the above inequality is minimized, that is, we fix $\mu$ to 
\[\mu_3^*:=\frac{2\lno \A\rno_{\R^n,\R^m}}{T+1}\sqrt{\frac{(1+\alpha)D_1}{D_2}},\]
and end up with the following bound:
\[\overline{\phi}(\bar{x})-\underline{\phi}(\bar{y})\leq \frac{4\lno\A\rno_{\R^n,\R^m}\sqrt{(1+\alpha)D_1D_2}}{T+1}+\frac{4(1+\alpha)D_1M}{(T+1)^2}.\]
Note that Alternative 2 coincides with Alternative 1 if $\alpha = 5(T+1)-1$. 

\section{An application in large-scale eigenvalue optimization}

In this section, we study the practical behavior of accelerated smoothing techniques. We apply them to the problem of finding a convex combination of given symmetric matrices such that the maximal eigenvalue of the resulting matrix is minimal. 

\subsection{Problem description}

Let \[\Delta_m:=\left\{ x\in\R^m_{\geq 0}:\sum_{j=1}^m x_j=1\right\}\subset\R^m\] be the $(m-1)$-dimensional probability simplex. Denoting by $\mathcal{S}_n$ the space of symmetric real $(n\times n)$-matrices, we write $Y\succeq 0$ if $Y\in \mathcal{S}_n$ is positive semidefinite and $\textup{Tr}(Y):=\sum_{i=1}^n Y_{ii}$ for the trace of $Y$. We refer to 
\[\Delta_n^M:=\left\{ Y\succeq 0: \textup{Tr}(Y)=1\right\}\subset\mathcal{S}_n\]  
as the simplex in matrix form. Finally, we denote by 
\[\lambda_n(Y)\geq \ldots\geq \lambda_1(Y)\]
the eigenvalues of the symmetric matrix $Y\succeq 0$ and assume that they are ordered decreasingly. Throughout this section, we consider the following problem:
\begin{eqnarray}\label{eq:min_max_eig}
\min_{x\in\Delta_m}\lambda_n\left(\sum_{j=1}^m x_jA_j\right)= \min_{x\in\Delta_m}\left\{\overline{\phi}(x):=\max_{Y\in\Delta_n^M}\sum_{j=1}^m x_j\la A_j, Y\ra_F\right\},  
\end{eqnarray}
where $A_1,\ldots,A_m\in\mathcal{S}_n$ and $\la\cdot,\cdot\ra_F$ denotes the Frobenius scalar product.

\subsection{Applying accelerated smoothing techniques}

\subsubsection{Smoothing the objective function}

We equip $\mathcal{S}_n$ with the induced $1$-norm, that is, with $\lno Y\rno_{(1)}:=\sum_{i=1}^n|\lambda_i(Y)|$, where $Y\in\mathcal{S}_n$. The dual norm corresponds to the induced $\infty$-norm, that is, to the norm $\lno W\rno_{(\infty)}:=\max_{1\leq i\leq n}|\lambda_i(W)|$ with $W\in\mathcal{S}_n$. We choose 
\[\dem (Y):=\ln(n)+\sum_{i=1}^n\lambda_i(Y)\ln(\lambda_i(Y)),\qquad Y\in\Delta_n^M,\]
as \dgf\ for the set $\Delta_n^M$, for which we have $\dem(Y)\leq \ln(n)$ for any $Y\in\Delta_n^M$; see for instance \cite{nesterov:coreDP73/2004} for a proof that $\dem$ is a \dgf\ for $\Delta_n^M$. We obtain the following smooth objective function as an approximation to $\overline{\phi}$:
\[\overline{\phi}_\mu(x):=\max_{Y\in\Delta_n^M}\left\{\sum_{j=1}^m x_j\la A_j, Y\ra_F-\mu\dem (Y) \right\}=\mu\ln\left(\sum_{i=1}^n\exp\left[\lambda_i\left(\frac{\sum_{j=1}^mx_jA_j}{\mu}\right)\right]\right)-\mu \ln(n),\]
where $x\in\Delta_m$ and $\mu>0$ denotes the smoothness parameter. The approximation quality depends on the smoothness parameter:
\[\overline{\phi}_\mu(x)\leq \overline{\phi}(x)\leq \overline{\phi}_\mu(x) +\mu\ln(n)\qquad\forall\ x\in\Delta_m.\]
Finally, the gradient of $\overline{\phi}_\mu$ is given by 
\[\left[\nabla \overline{\phi}_\mu(x)\right]_j=\la A_j, Y_*(x)\ra_F\qquad\forall\ 1\leq j\leq m,\]
where $x\in\Delta_m$ and $Y_*(x)$ denotes the unique maximizer of $Y\mapsto \sum_{j=1}^n x_j\la A_j, Y\ra_F-\mu \dem(Y)$ over $\Delta_n^M$. Theorem \ref{thm:gradient} implies that the gradient is Lipschitz continuous with a Lipschitz constant of $L_\mu:=\max_{1\leq j\leq m}\lno A_j\rno_{(\infty)}/\mu$.

\subsubsection{Applying the accelerated optimal first-order method with hybrid setting}

Let the space $\R^m$ be equipped with the $1$-norm. We use 
\[\de (x):=\ln(m)+\sum_{j=1}^mx_j\ln(x_j),\qquad x\in\Delta_m,\]
as \dgf\ for the set $\Delta_m$. Note that $\de(x)\leq \ln(m)$ for any $x\in\Delta_m$. 

We run Algorithm \ref{alg:smoothing_techniques} with the hybrid setting that is described in Alternative 2 in Section \ref{sec:ST}. Let us fix the accuracy $\epsilon>0$ and the parameter $\alpha\geq 0$ that defines when to switch back to the non-adaptive setting. The smoothness parameter is set as follows:
\[\mu:=\frac{\epsilon}{2\ln(n)}.\]
Note that the smoothness parameter does not depend on $\alpha$. According to Theorem \ref{thm:convergence_smoothing}, we need to perform at most
\begin{eqnarray}\label{eq:worst_case_T}
T=\left\lceil\frac{4\max_{1\leq j\leq m}\lno A_j\rno_{(\infty)}\sqrt{(1+\alpha)\ln(m)\ln(n)}}{\epsilon}-1\right\rceil 
\end{eqnarray}
iterations of Algorithm \ref{alg:smoothing_techniques} in order to find a tuple $(\bar{x},\bar{Y})\in\Delta_m\times\Delta_n^M$ such that
\begin{eqnarray}\label{eq:duality_gap}\max_{Y\in\Delta_n^M}\sum_{j=1}^m \bar x_j\la A_j, Y\ra_F - \min_{x\in\Delta_m}\sum_{j=1}^m x_j\la A_j,\bar Y\ra_F\leq \epsilon.\end{eqnarray}

\subsection{Numerical results}

We consider randomly generated instances of Problem (\ref{eq:problem_min_max}), where we fix $m$ to $100$ and where the symmetric $(n\times n)$-matrices $A_1,\ldots,A_m$ have a joint sparsity structure, each of them with about $n^2/10$ non-zero entries. We approximate the parameter \[\L:=\max_{1\leq j\leq m}\lno A_j\rno_{(\infty)}\] by applying the Power method to the matrices $A_j$ and taking the maximum, which we denote by $\L'$, of the computed values afterwards. We solve the randomly generated instances of Problem (\ref{eq:problem_min_max}) up to a relative accuracy of $\epsilon\L'$ with $\epsilon:=0.002$.

All numerical results that we present in this section are averaged over ten runs and obtained on a computer with 24 processors, each of them with 2.67 GHz, and with 96 GB of RAM. The methods are implemented in Matlab (version R2012a). Matrix exponentials are computed through the Matlab built-in function \texttt{expm()}.

\subsubsection{Comparing the practical behavior of different methods}

In Table \ref{table:different_methods}, we present numerical results for the following two methods:
\begin{itemize}
 \item Original smoothing techniques: This implementation corresponds to Algorithm \ref{alg:smoothing_techniques} with constant $L_t=L_\mu$ for any $0\leq t\leq T$. That is, we set $\alpha=0$ in Alternative 2 in Section \ref{sec:ST}.
 \item Accelerated smoothing techniques: We equip Algorithm \ref{alg:smoothing_techniques} with the hybrid setting described in Alternative 2 in Section \ref{sec:ST}, where we choose $\alpha:=3$ and $\kappa:=10^{-12}$. With this setting, we need to perform twice as many iterations as with original smoothing techniques with respect to the worst-case bounds; see (\ref{eq:worst_case_T}).
\end{itemize}
For both methods, we check the duality gap (\ref{eq:duality_gap}) at every $100$-th iteration. Additionally for the later method, we also verify this condition at every of the first hundred iterations. The maximal eigenvalue that corresponds to the first term in (\ref{eq:duality_gap}) is computed through the Matlab built-in functions $\texttt{max()}$ and $\texttt{eig()}$.

\begin{table}[t]
\begin{center}
\textbf{Average CPU time [sec]}\\
\begin{tabular}{|l||cccc|}
 \hline 
$n$ & $100$ & $200$ & $400$ &$800$ \cr
\hline
Original smoothing techniques & $139$ & $366$ & $1'406$ & $5'961$ \cr
Accelerated smoothing techniques & $116$ & $3$ & $9$ & $32$ \cr
\hline
Acceleration & $16.55\%$ & $99.18\%$ & $99.36\%$ & $99.46\%$\cr
\hline
 \end{tabular}\\
\vspace{2mm}
\textbf{Average \# of iterations that are required in practice}\\
\begin{tabular}{|l||cccc|}
 \hline 
$n$ & $100$ & $200$ & $400$ &$800$  \cr
\hline
Original smoothing techniques & $6'180$ &$6'690$ & $7'150$ & $7'520$  \cr
Accelerated smoothing techniques & $4'918$ & $18$ & $14$ & $13$\cr
\hline
Reduction & $20.42\%$ & $99.73\%$ & $99.80\%$ & $99.83\%$ \cr
\hline
 \end{tabular}\\
\vspace{2mm}
\textbf{Average \# of iterations that are required in theory}\\
\begin{tabular}{|l||cccc|}
 \hline 
$n$  & $100$ & $200$ & $400$ &$800$  \cr
\hline
Original smoothing techniques& $9'210$ & $9'879$  & $10'505$ & $11'096$ \cr
Accelerated smoothing techniques& $18'420$ & $19'758$ & $21'011$ & $22'193$\cr
\hline
Reduction & $-100.00\%$ & $-100.00\%$ & $-100.01\%$ & $-100.01\%$\cr
\hline
 \end{tabular}
\end{center}
\caption{Average CPU time and number of iterations (in practice and in theory) that are required by original and accelerated smoothing techniques for finding an approximate solution to randomly generated instances of Problem (\ref{eq:problem_min_max}) (with fixed accuracy $0.002\L'$ and with $m=100$).}
\label{table:different_methods}
\end{table}

We observe that accelerated smoothing techniques require significantly less CPU time and iterations in practice than original smoothing techniques; see Table \ref{table:different_methods}. For problems involving matrices of size $200\times 200$ up to size $800\times 800$, we can reduce the number of iterations in practice and the CPU time by more than $99\%$. Interestingly, the number of iterations that are required by accelerated smoothing techniques in practice is even decaying when the matrix size $n$ is getting larger.

\begin{figure}[h]
\begin{center}
\includegraphics[width=0.48\linewidth]{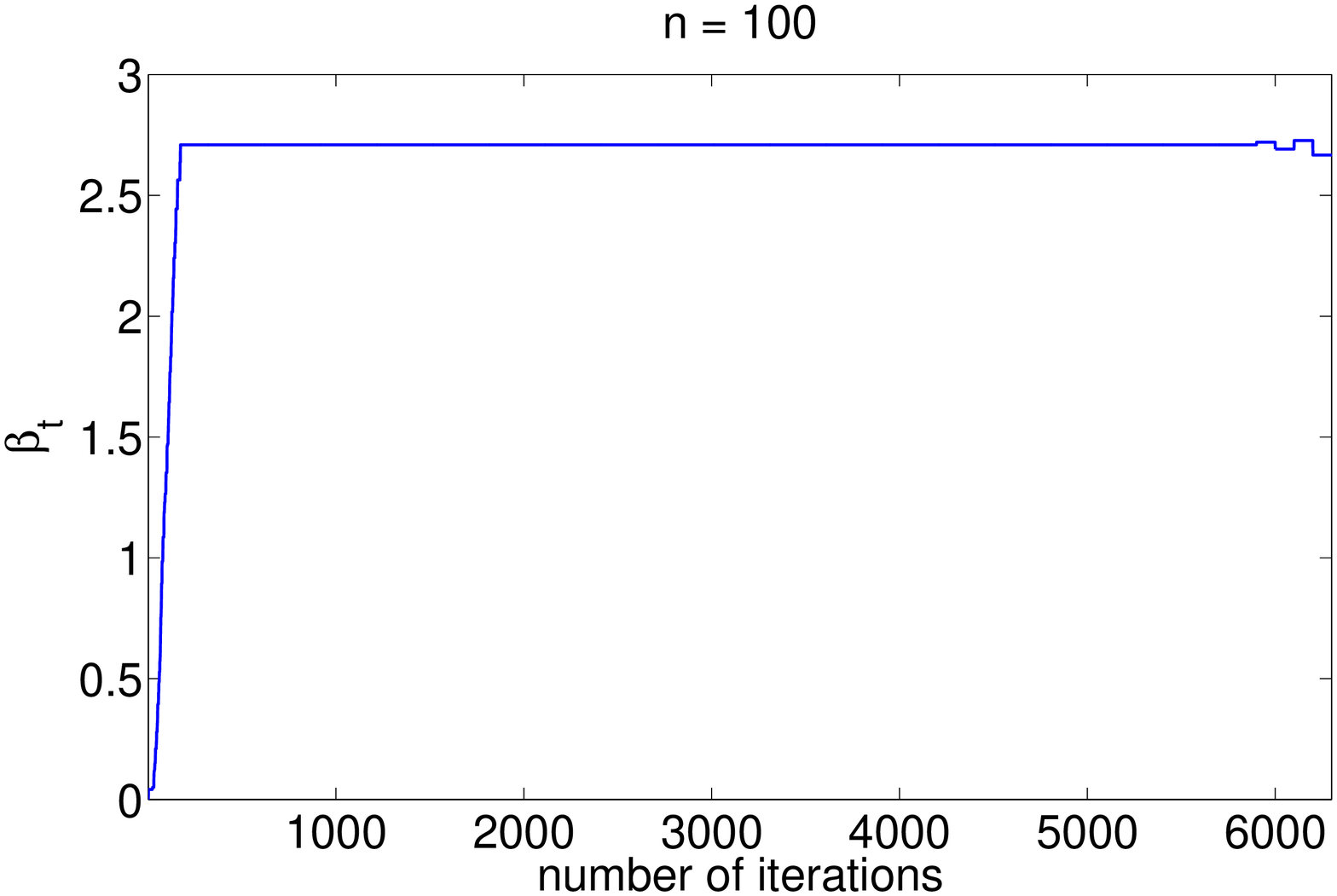}
\includegraphics[width=0.48\linewidth]{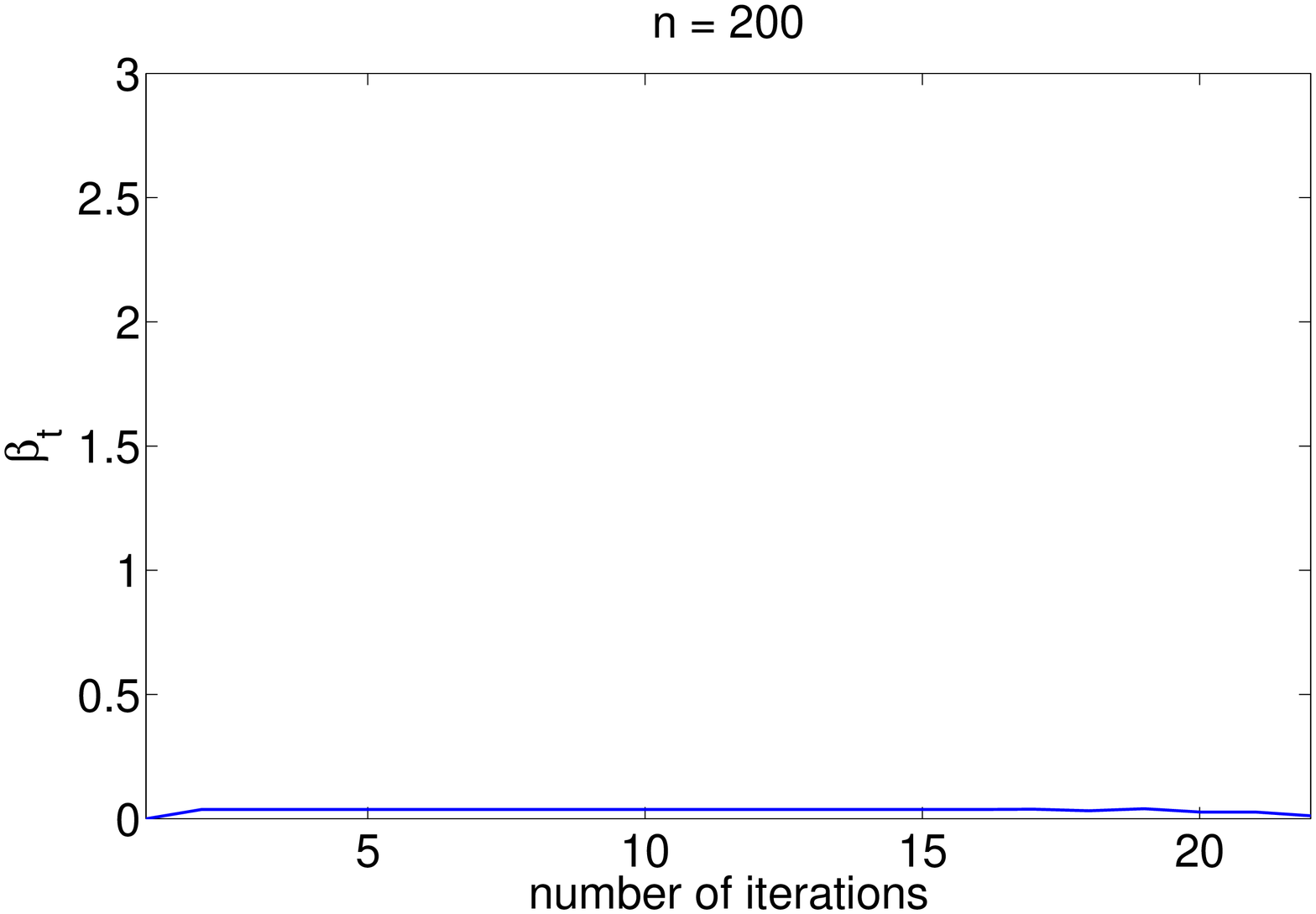}
\includegraphics[width=0.48\linewidth]{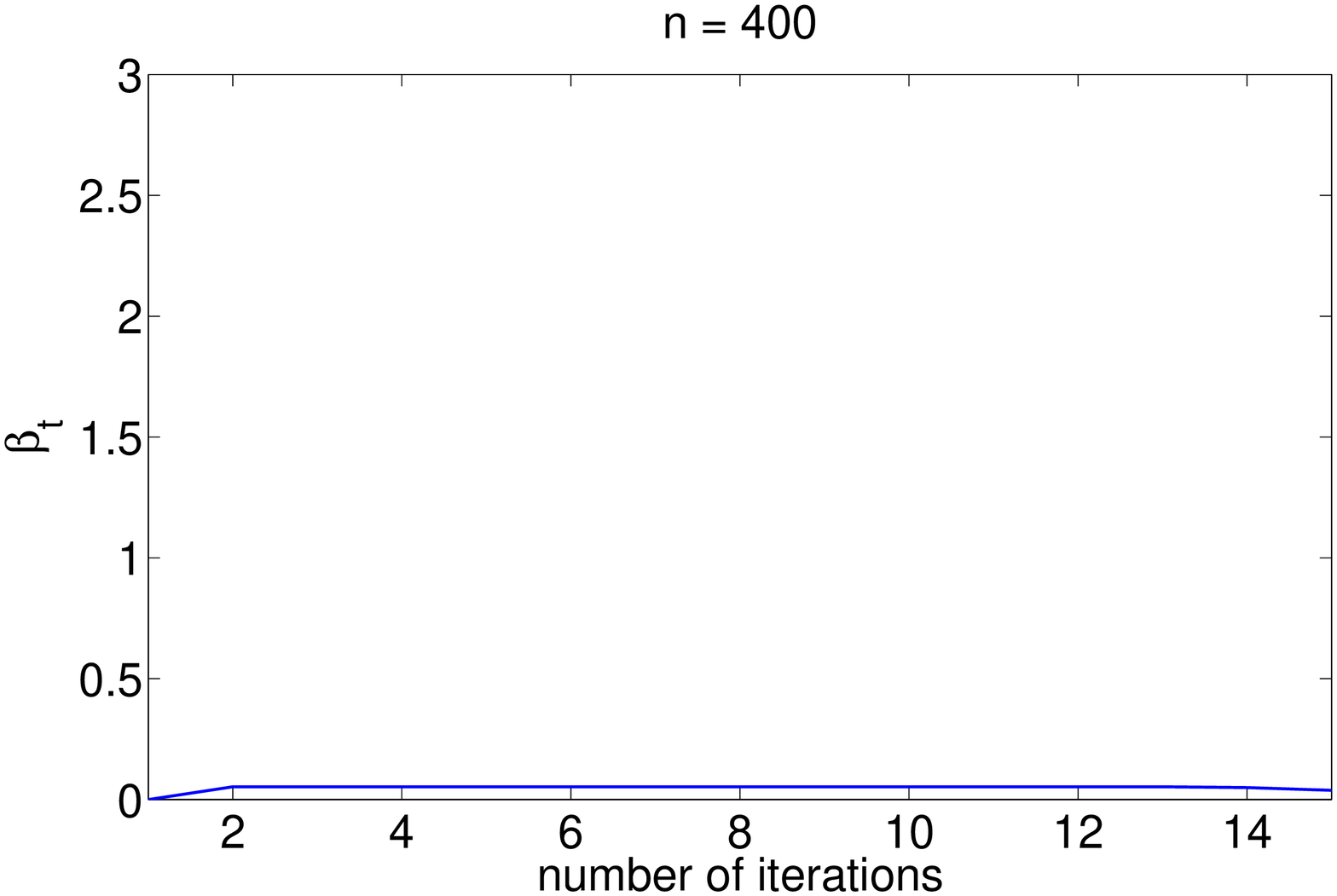}
\includegraphics[width=0.48\linewidth]{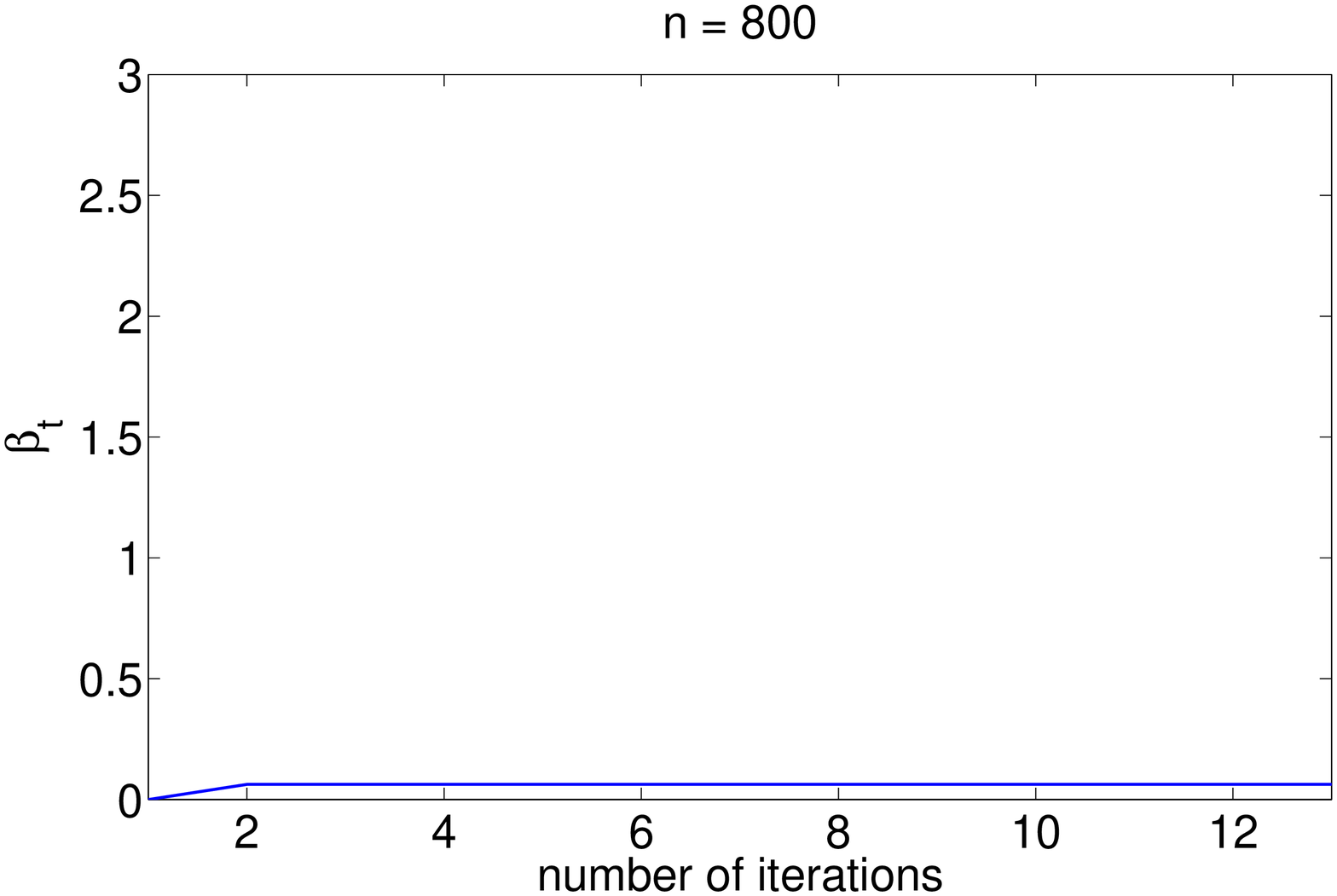}
\end{center}
\caption{Ratios $\beta_t$; see (\ref{eq:beta}) for the definition of these ratios.}
\label{fig:beta}
\end{figure}

Note that there exists a gap in the average CPU time and number of iterations that are required by accelerated smoothing techniques in practice for solving the instances of size $100\times 100$ and the instances of size $200\times 200$. In Figure \ref{fig:beta}, we plot the values 
\begin{eqnarray}\label{eq:beta}\beta_t:=\frac{-\chi_t}{\ln(m)L_0}= \frac{-\sum_{t'=0}^{t-1}\left(L_{t'+1}-L_{t'}\right)\left(\de (z_{t'+1})-\frac{1}{2}\left\|z_{t'}-\hat{x}_{t'+1}\right\|_{1}^2 \right)}{D_1L_0}\qquad \forall\ t\geq 1.\end{eqnarray}
In contrast to the cases $n=200$, $n=400$, or $n=800$, where these values remain small (that is, below $0.25$), we have considerably large values $\beta_t$ for $n=100$. However, the values are still below $3$, as we switch back to a non-adaptive setting as soon as $\beta_t$ would be larger than $3$. This behavior is in full accordance with the gap mentioned in the beginning of this paragraph. The non-smooth patterns at the end of the plots in Figure \ref{fig:beta} are due to the averaging over the different runs (We may need a different number of iterations in the different runs.).

\subsubsection{Solving problems of very large scale}

In Table \ref{table:large_scale}, we show numerical results for accelerated smoothing techniques (with $\alpha=3$, $\kappa=10^{-12}$, and the same duality gap checking procedure as above) when applied to randomly generated instances of (\ref{eq:min_max_eig}) that are of very large scale. Using accelerated smoothing techniques, we are able to solve approximately instances of (\ref{eq:min_max_eig}) involving matrices of size $12'800\times 12'800$ in about $8$ hours and $40$ minutes on average. Clearly, this performance would be out of reach for original smoothing techniques.

\begin{table}[t]
\begin{center}
\textbf{Accelerated smoothing techniques applied to large-scale instances of (\ref{eq:min_max_eig})}\\
\begin{tabular}{|l||cccc|}
 \hline 
$n$ & $1'600$ & $3'200$ & $6'400$ &$12'800$ \cr
\hline
CPU time [sec] & $158$ & $791$ & $4'566$ & $31'240$ \cr
Average \# of iterations that are required in practice & $13$ & $13$ & $13$ & $13$ \cr
Average \# of iterations that are required in theory & $23'315$ & $24'386$ & $25'411$ & $26'397$ \cr
\hline
 \end{tabular}\\
\end{center}
\caption{Average CPU time and number of iterations (in practice and in theory) that are required by accelerated smoothing techniques for finding an approximate solution to randomly generated large-scale instances of Problem (\ref{eq:problem_min_max}) (with fixed accuracy $0.002\L'$ and with $m=100$).}
\label{table:large_scale}
\end{table}

\section*{Acknowledgments}
We gratefully thank Yurii Nesterov and Hans-Jakob L\"uthi for many helpful discussions. This research is partially funded by the Swiss National Fund.

\appendix

\section{Proof of Theorem \ref{thm:I_t}}
\label{sec:proof_opt_first_order}

Choose $T\in\N_0$ and let the sequences $\left(x_t\right)_{t=0}^{T+1}$, $\left(u_t\right)_{t=0}^{T+1}$, $\left(z_t\right)_{t=0}^T$, $\left(\hat{x}_t\right)_{t= 1}^{T+1}$, $\left(\gamma_t\right)_{t=0}^{T+1}$, $\left(\sumweights_t\right)_{t= 0}^{T+1}$, $\left(\tau_t\right)_{t=0}^T$, and $\left(L_t\right)_{t=0}^T$ be generated by Algorithm \ref{alg:opt_first_order}. Recall that Inequality $(\mathcal{I}_t)$ holds for $0\leq t\leq T$ if 
\begin{equation}
\sumweights_tf(u_t)+\sum_{k=0}^{t-1}\left(L_{k+1}-L_k\right)\left(d_Q(z_{k+1})-\frac{1}{2}\left\|z_k-\hat{x}_{k+1}\right\|_{\R^n}^2 \right)\leq \psi_t, \tag{$\mathcal{I}_t$}
\end{equation}
where
\[\psi_t:=\min_{x\in Q}\left\lbrace
       \sum_{k=0}^{t}\gamma_k\left(f(x_k)+\langle \nabla f(x_k), x-x_k\rangle \right) + L_td_Q(x) \right\rbrace.\]
By its definition (see Algorithm \ref{alg:opt_first_order}), the element $z_t\in Q$ is the minimizer to the above optimization problem, which allows us to rewrite $\psi_t$ as:
\[\psi_t=\sum_{k=0}^{t}\gamma_k\left(f(x_k)+\langle \nabla f(x_k), z_t-x_k\rangle \right)+L_td_Q(z_t).\]

We show by induction that Inequality $(\mathcal{I}_t)$ holds for any $0\leq t\leq T$.

\begin{lem} \label{lem:basic_step}
Inequality $(\mathcal{I}_0)$ holds, that is, we have $\gamma_0f(u_0)\leq \psi_0$.
\end{lem}

\noindent\textbf{Proof:} We apply the definition of $u_0$ (see Algorithm \ref{alg:opt_first_order}), Inequality (\ref{eq:lower_bound_dgf}), the condition on $\gamma_0$ saying that $\gamma_0\in(0,1]$, and Theorem 2.1.5 in \cite{nesterovintrodlectures} in order to justify the following relations:
\begin{eqnarray*}
 \psi_0&:=&\min_{x\in Q}\left\lbrace
        \gamma_0\left(f(x_0)+\langle \nabla f(x_0, x-x_0\rangle \right) +L_0d_Q(x)\right\rbrace\cr
  &=&\gamma_0\left(f(x_0)+\langle \nabla f(x_0), u_0-x_0\rangle\right) +L_0d_Q(u_0) \cr
 &\geq&\gamma_0\left(f(x_0)+\langle \nabla f(x_0), u_0-x_0\rangle \right) +\frac{L_0}{2}\left\|u_0-x_0\right\|_{\R^n}^2 \cr
&\geq&\gamma_0\left( f(x_0)+\langle \nabla f(x_0), u_0-x_0\rangle+\frac{L_0}{2}\left\|u_0-x_0\right\|_{\R^n}^2  \right)\cr
&\geq& \gamma_0f(u_0).
\end{eqnarray*}
\qed

Let us verify the inductive step.

\begin{lem} \label{lem:inductive_step}
Let $0\leq t\leq T-1$. If Inequality $(\mathcal{I}_t)$ holds, also $(\mathcal{I}_{t+1})$ is true. 
\end{lem}

\noindent\textbf{Proof:} Let $0\leq t\leq T-1$ and assume that $(\mathcal{I}_t)$ holds. We make the following two definitions:
\begin{eqnarray*}
\chi_t&:=&\sum_{k=0}^{t-1}\left(L_{k+1}-L_k\right)\left(d_Q(z_{k+1})-\frac{1}{2}\left\|z_k-\hat{x}_{k+1}\right\|_{\R^n}^2 \right)\in\R,\cr
s_{t}&:=&\sum_{k=0}^t\gamma_k\nabla f(x_k)\in \R^n.
\end{eqnarray*}
In addition, we define the linear function:
\begin{eqnarray*}
l_t:Q\rightarrow \R:x\mapsto l_t(x)=\sum_{k=0}^t\gamma_k\left( f(x_k)+\left\langle \nabla f(x_k),x-x_k\right\rangle\right).
\end{eqnarray*}
Choose $x\in Q$. The definition of $z_t$ implies:
\begin{eqnarray}
 \label{eq:z_k_minimizer}
0\leq \left\langle L_t\nabla d_Q(z_t)+\sum_{k=0}^t \gamma_k\nabla f(x_k), x-z_t\right\rangle= \left\langle L_t\nabla d_Q(z_t)+s_t, x-z_t\right\rangle.
\end{eqnarray}
As the Inequality $(\mathcal{I}_t)$ holds and as the function $f$ is convex, we have:
\begin{eqnarray*}
\psi_t &\geq& \sumweights_t f(u_t)+\chi_t \geq \sumweights_t \left(f(x_{t+1})+\left\langle \nabla f(x_{t+1}),u_t-x_{t+1}\right\rangle \right)+\chi_t.
\end{eqnarray*}
This implies:
\begin{eqnarray*}
\hspace{-1cm}\psi_t +\gamma_{t+1}\left(f(x_{t+1})+\left\langle \nabla f(x_{t+1}),x-x_{t+1}\right\rangle\right) \geq \sumweights_{t+1}f(x_{t+1})+\gamma_{t+1}\left\langle \nabla f(x_{t+1}),x-z_t\right\rangle+\chi_t,
\end{eqnarray*}
where we use the relations $\sumweights_{t+1}=\sumweights_t+\gamma_{t+1}$ and
\begin{eqnarray*}
\sumweights_t(u_t-x_{t+1})+\gamma_{t+1}(x-x_{t+1})&=&\sumweights_tu_t-\sumweights_{t+1}x_{t+1}+\gamma_{t+1}x\cr
&=&\sumweights_tu_t-\sumweights_{t+1}\left(\tau_tz_t+(1-\tau_t)u_t\right)+\gamma_{t+1}x\cr
&=&\sumweights_tu_t-\sumweights_{t+1}\left(\frac{\gamma_{t+1}}{\sumweights_{t+1}}z_t+\frac{\sumweights_t}{\sumweights_{t+1}}u_t\right)+\gamma_{t+1}x\cr
&=&\gamma_{t+1}(x-z_t).
\end{eqnarray*}
Combining the above inequality with the fact that $\psi_t=L_td_Q(z_t)+l_t(z_t)$ and with (\ref{eq:z_k_minimizer}), we observe:
\begin{eqnarray*}
&&\hspace{-1cm} L_td_Q(x)+l_{t+1}(x)\cr
&=&L_td_Q(x)+l_t(x)+\gamma_{t+1}\left(f(x_{t+1})+\left\langle \nabla f(x_{t+1}),x-x_{t+1}\right\rangle\right)\cr
&=&L_t V_{z_t}(x)+\psi_t+\left\langle L_t\nabla d_Q(z_t)+s_t,x-z_t\right\rangle+\gamma_{t+1}\left(\left\langle \nabla f(x_{t+1}),x-x_{t+1}\right\rangle+f(x_{t+1})\right)\cr
&\geq& L_t V_{z_t}(x)+\psi_t+\gamma_{t+1}\left(f(x_{t+1})+\left\langle \nabla f(x_{t+1}),x-x_{t+1}\right\rangle\right)\cr
&\geq& L_t V_{z_t}(x)+\sumweights_{t+1}f(x_{t+1})+\gamma_{t+1}\left\langle \nabla f(x_{t+1}),x-z_t\right\rangle+\chi_t.
\end{eqnarray*}
With $\vartheta_t^{(1)}:= \left(L_{t+1}-L_t\right)d_Q(z_{t+1})$, we thus get:
\begin{eqnarray*}
\psi_{t+1}&:=&\min_{x\in Q}\left\lbrace
        L_{t+1}d_Q(x)+l_{t+1}(x) \right\rbrace\cr
&=&L_{t+1}d_Q(z_{t+1})+l_{t+1}(z_{t+1})\cr
&=&\vartheta_t^{(1)}+L_td_Q(z_{t+1})+l_{t+1}(z_{t+1})\cr
&\geq&\vartheta_t^{(1)}+\min_{x\in Q}\left\{L_td_Q(x)+l_{t+1}(x)\right\}\cr
&\geq&\vartheta_t^{(1)}+\min_{x\in Q}\left\{L_t V_{z_t}(x)+\sumweights_{t+1}f(x_{t+1})+\gamma_{t+1}\left\langle \nabla f(x_{t+1}),x-z_t\right\rangle+\chi_t\right\}.
\end{eqnarray*}
Let $\vartheta_t^{(2)}:=\frac{1}{2}\left(L_t-L_{t+1}\right)\left\|z_t-\hat{x}_{t+1}\right\|_{\R^n}^2$. Using the construction rule for $\hat{x}_{t+1}$ and the fact that the inequality $V_z(x)\geq \lno x-z\rno_{\R^n}^2/2$ holds for any $x\in Q$ and $z\in Q^o$ (this relation follows from the strong convexity of $d_Q$), we obtain:
\begin{eqnarray*}
\psi_{t+1}&\geq&\vartheta_t^{(1)}+L_t V_{z_t}(\hat{x}_{t+1})+\sumweights_{t+1}f(x_{t+1})+\gamma_{t+1}\left\langle \nabla f(x_{t+1}),\hat{x}_{t+1}-z_t\right\rangle+\chi_t\cr
&\geq&\vartheta_t^{(1)}+\frac{L_t}{2}\left\|z_t-\hat{x}_{t+1}\right\|_{\R^n}^2+\sumweights_{t+1}f(x_{t+1})+\gamma_{t+1}\left\langle \nabla f(x_{t+1}),\hat{x}_{t+1}-z_t\right\rangle+\chi_t\cr
&=&\vartheta_t^{(1)}+\vartheta_t^{(2)} + \frac{L_{t+1}}{2}\left\|z_t-\hat{x}_{t+1}\right\|_{\R^n}^2+\sumweights_{t+1}f(x_{t+1})+\gamma_{t+1}\left\langle \nabla f(x_{t+1}),\hat{x}_{t+1}-z_{t}\right\rangle+\chi_t\cr
&=&\vartheta_t^{(1)}+\vartheta_t^{(2)}+\chi_t + \sumweights_{t+1}\left( \frac{L_{t+1}}{2\sumweights_{t+1}}\left\|z_t-\hat{x}_{t+1}\right\|_{\R^n}^2+f(x_{t+1})+\tau_t\left\langle \nabla f(x_{t+1}),\hat{x}_{t+1}-z_t\right\rangle\right).
\end{eqnarray*}
As $\tau_t^2\leq \sumweights_{t+1}^{-1}$ and as $x_{t+1}-\tau_tz_t=(1-\tau_t)u_t=u_{t+1}-\tau_t\hat{x}_{t+1}$, this inequality yields to:
\begin{eqnarray*}
\psi_{t+1}&\geq&\vartheta_t^{(1)}+\vartheta_t^{(2)} +\chi_t + \sumweights_{t+1}\left( \frac{L_{t+1}\tau_t^2}{2}\left\|z_t-\hat{x}_{t+1}\right\|_{\R^n}^2+f(x_{t+1})+\tau_t\left\langle \nabla f(x_{t+1}),\hat{x}_{t+1}-z_t\right\rangle\right)\cr
&=&\vartheta_t^{(1)}+\vartheta_t^{(2)}+\chi_t + \sumweights_{t+1}\left( \frac{L_{t+1}}{2}\left\|u_{t+1}-x_{t+1}\right\|_{\R^n}^2+f(x_{t+1})+\left\langle \nabla f(x_{t+1}),u_{t+1}-x_{t+1}\right\rangle\right).
\end{eqnarray*}
It remains to apply (\ref{eq:adapted_L_cond_opt}):
\begin{eqnarray*}
 \psi_{t+1}&\geq&\vartheta_t^{(1)}+\vartheta_t^{(2)} + \sumweights_{t+1}f(u_{t+1})+\chi_t\cr
&=&\sum_{k=0}^{t}\left(L_{k+1}-L_k\right)\left(d_Q(z_{k+1})-\frac{1}{2}\left\|z_k-\hat{x}_{k+1}\right\|_{\R^n}^2  \right)+ \sumweights_{t+1}f(u_{t+1}).
\end{eqnarray*}
\qed

\bibliographystyle{amsalpha}
\bibliography{Jordan}

\end{document}